\documentclass{article}

\usepackage{amsfonts}
\usepackage{amsmath}
\usepackage{amssymb}
\usepackage{float} 
\usepackage{subcaption}		
\usepackage{graphicx}
\usepackage{verbatim}
\usepackage[Algorithm,ruled]{algorithm}       
\usepackage{algpseudocode}
\usepackage{pbox}
\usepackage{tabularx}
\usepackage{placeins} 
\usepackage[margin =1.0in]{geometry}

\newcommand{\Om}{\Omega}
\newcommand{\DOm}{\partial\Omega}

\newcommand{\R}{\mathbb{R}}
\renewcommand{\Re}{\mathrm{Re}}

\newcommand{\conj}{\overline}
\newcommand{\T}{{\mathbf{t}}}

\newcommand{\Texp}{\mathbf{t}^{\mbox{{\tiny \rm exp}}}}

\newcommand{\Tpr}{\T_{\rm pr}}
\newcommand{\TRBIE}{\T_{\mbox{\tiny $R$}}^{\mbox{{\tiny \rm BIE}}}}

\newcommand{\TRs}{\T_{\mbox{\tiny $R_1$,$R_2$}}}
\newcommand{\Lam}{\Lambda}
\newcommand{\Lampr}{\Lambda_{\rm pr}}
\newcommand{\dbar}{\bar{\partial}}

\newcommand{\sigpr}{\sigma_{\rm pr}}
\newcommand{\tsigpr}{\tilde{\sigma}_{\rm pr}}

\newcommand{\sigalpha}{\sigma_{\mbox{\tiny $R_2$,$\alpha$}}}
\newcommand{\mupr}{\mu_{\rm pr}}
\newcommand{\muint}{\mu_{\rm int}}
\newcommand{\mutilde}{\tilde{\mu}_{\mbox{\tiny $R_2$}}}

\newcommand{\mualpha}{\mu_{\mbox{\tiny $R_2$,$\alpha$}}}
\newcommand{\psipr}{\psi_{\rm pr}}

\newcommand{\qpr}{q_{\rm pr}}

\algsetblock{Proc}{EndProc}{100}{0.8cm}
\algrenewtext{Proc}[1]{{#1}}
\algnotext{EndProc}
\newcommand{\DashState}{\State- }

\title{A D-BAR ALGORITHM WITH A PRIORI INFORMATION FOR 2-D ELECTRICAL IMPEDANCE TOMOGRAPHY} 

\author{Melody Alsaker\thanks{Department of Mathematics, Colorado State University, USA} \and Jennifer L. Mueller\thanks{Department of Mathematics and School of Biomedical Engineering, Colorado State University, USA}}

\begin{document}
\maketitle

\begin{abstract}
A method for including {\em a priori} information in the 2-D D-bar algorithm is presented.  Two methods of assigning conductivity values to the prior are presented, each corresponding to a different scenario on applications.  The method is tested on several numerical examples with and without noise and is demonstrated to be highly effective in improving the spatial resolution of the D-bar method.
\end{abstract}

\pagestyle{myheadings}
\thispagestyle{plain}
\markboth{M. ALSAKER AND J. L. MUELLER}{D-BAR WITH A PRIORI INFORMATION FOR EIT}

\section{Introduction}

Electrical impedance tomography (EIT) is a low-cost, portable, and noninvasive imaging modality that is free of ionizing radiation with many potential applications for pulmonary imaging. In EIT, an image is formed of the conductivity distribution $\sigma$ inside a body using only surface voltage and current measurements.  Mathematically, this is a nonlinear inverse problem which is well known to be extremely ill-posed.  A significant challenge in EIT imaging is the computation of static images with high-quality spatial resolution. Due to the ill-posedness, finer details in the image are often lost in the presence of noisy measurements.  Including prior information in the reconstruction algorithm has been shown to be one way to improve spatial resolution \cite{Avis95, Baysal98, CamargoThesis, Dehghani99, Dobson94,  Ferrario12, FloresTapia10, Soleimani06, Vauhkonen98}.  This prior knowledge corresponds to a clinical situation in which we have a CT scan (or other similar data) for a human subject from which we may extract information regarding spatial locations of organ boundaries or conductivity estimates. When diagnosing and treating certain lung conditions, it is often necessary to obtain repeated thoracic CT scans, each of which imparts a dose of harmful radiation. EIT scans, on the other hand, have no ill effects. It is therefore highly desirable to use \emph{a priori} information obtained from a CT or other scan to provide an improved EIT image, and then perform repeated harmless and comparatively inexpensive EIT scans in place of follow-up CT scans.  

Reconstruction algorithms that involve the minimization of a cost functional, such as a Gauss-Newton algorithm, include the {\em a priori} information in the penalty term, penalizing reconstructions that deviate too greatly from the prior in a given norm. This technique does not generalize to noniterative algorithms, and until now there has been no direct reconstruction method to utilize \emph{a priori} information. The algorithm presented here therefore represents the first direct reconstruction method for EIT to incorporate \emph{a priori} data.

In this paper, we first provide an outline of the D-bar method without {\em a priori} data  in \S\ref{Background}. The \emph{a priori} scheme for D-bar methods is described in \S\ref{APriori_outline}, in which spatial information regarding locations of inclusion boundaries and approximated conductivity values are encoded into the equations for D-bar.  Results using simulated data with and without noise on a circular domain with adjacent current patterns and 32 electrodes are presented in \S\ref{Results}.

\section{Theoretical background} \label{Background}
%

In EIT,  current is applied on electrodes on the surface of a domain, and  the resulting  voltages are measured on the electrodes.  Let $\Om \subset \R^2$ be a bounded, simply connected Lipschitz domain. Then the electric potential $u(x,y)$ within $\Om$ is modeled by
\begin{equation}\label{cond_eqn}
\nabla \cdot (\sigma(x,y) \nabla u(x,y) )= 0, \quad (x,y) \in  \Om,
\end{equation} 
where $\sigma(x,y)$ is the conductivity distribution within $\Om$. The boundary data for the inverse problem is given by the Dirichlet-to-Neumann (DN) map $\Lam_\sigma$, which takes the boundary voltages to the current densities on the boundary:
\begin{equation}\label{DN_map}
\Lam_\sigma : u(x,y)\big|_{\DOm} \to \sigma(x,y) \frac{\partial u}{\partial \nu} \bigg|_{\DOm}, 
\end{equation}
where $\nu$ is the outward normal to the surface. The inverse conductivity problem is to reconstruct $\sigma(x,y)$ for $(x,y) \in \Om$ given knowledge of the DN map $\Lam_\sigma$.   We will denote $\Lam_1$ to be the DN map corresponding to the case of constant conductivity $\sigma \equiv 1$.

In this work, we will modify the 2-D D-bar reconstruction method based on the constructive global uniqueness proof in \cite{Nachman96} to include {\em a priori} information about the conductivity.  In \cite{Nachman96} it is shown that twice-differentiable conductivities $\sigma$ can be uniquely determined from $\Lam_\sigma$.   It is well-known that for $\sigma \in C^2(\Om)$, equation \eqref{cond_eqn} can be transformed to the Schr\"odinger equation through the change of variables $q =  \frac{\Delta \sqrt{\sigma}}{\sqrt{\sigma}}$, $\tilde u = \sqrt{\sigma}.$  Under the assumption that $\sigma$ is constant in a neighborhood of $\DOm$, one can smoothly extend $q = 0$  to $\R^2$ and consider the Schr\"odinger equation in the entire plane.  Points $(x,y)\in \R^2$ will be identified with points $z=x+iy\in \mathbb{C}$.  The direct reconstruction method relies on exponentially growing solutions to this Schr\"odinger equation,  also known as complex geometrical optics (CGO) solutions, which arise when a complex frequency parameter $k = k_1 + i k_2 \in \mathbb C$ is introduced in the equation and one seeks  solutions $\psi(z,k)$ satisfying the asymptotic condition $ e^{-i k z}\psi(z,k)-1 \in W^{1,p}(\R^2), p>2$.  It was proved in \cite{Nachman96} that the Schr\"odinger equation
\begin{equation}\label{Schrod_eqn}
(-\Delta + q(z)) \psi(z,k) = 0, \quad z \in \R^2
\end{equation}
has a unique solution with the desired asymptotic property when $q$ is of the form $\frac{\Delta \sqrt{\sigma(z)}}{\sqrt{\sigma(z)}}$.

The CGO solution $\psi(z,k)$ and its relative  $\mu(z,k) := e^{-i k z}\psi(z,k)$ are key to the direct reconstruction algorithm.  The algorithm has been implemented numerically and tested on simulated and experimental data \cite{DeAngelo10,Dodd14,Isaacson04,Isaacson06, Mueller03,Knudsen07,Murphy09,Siltanen00}.  A nonlinear regularization method for this algorithm was provided with proof in \cite{Knudsen09}.  Since details of the algorithm appear in numerous places in the literature, such as the papers above and \cite{MuellerBook}, we will give only a brief summary here with a focus on computation and the regularized method.  As in \cite{Nachman96} we will assume without loss of generality in this section that $\sigma=1$ in a neighborhood of $\DOm$.

The steps of the algorithm are

\noindent
Step 1.  Compute the CGO solution $\psi$ on $\DOm$ from the DN data by solving the boundary integral equation
\begin{equation}\label{BIE}
\psi(z,k)|_{\DOm} = e^{ikz}|_{\DOm} - \int_{\DOm} G_k(z - \zeta)(\Lam_\sigma - \Lam_1)\psi(\cdot,k) ds,
\end{equation}
where $G_k$ is the Faddeev Green's function \cite{Faddeev66} for the Laplacian, defined by 
\begin{equation} \label{Gk}
G_k(z) := e^{ikz} \int_{\R^2} \frac{e^{iz \cdot \xi}}{\xi (\bar \xi + 2k)} d\xi, \quad -\Delta G_k(z) = \delta(z).
\end{equation}

\noindent
Step 2.  Compute the scattering transform from 
\begin{equation}\label{scatt_trans2}
\T(k) =  \int_{\DOm} e^{i \bar{k} \bar{z}} (\Lam_\sigma - \Lam_1) \psi(z,k) \; ds.
\end{equation}

\noindent
Step 3.  Solve the D-bar ($\dbar$) equation for $\mu(z,k)$
\begin{equation} \label{Dbar_eq}
\dbar_k \mu(z,k) = \frac{\T(k)}{4 \pi \bar k} e^{-i(kz + \bar k \bar z)} \conj{\mu(z,k)}, 
\end{equation}
 which can be written in integral form as
\begin{equation}\label{Dbar_int_eqn}
\mu(z,k) = 1 + \frac{1}{(2\pi)^2} \int_{\R^2} \frac{\T(k')}{\bar{k'}(k-k')} e^{-i (zk' + \bar{z}\bar{k'})} \overline{\mu(z,k')} \; dk'.
\end{equation}

\noindent
Step 4.  Compute the conductivity $\sigma(z)$ for each $z\in\Om$ or in a region of interest from
\begin{equation}\label{sigma_mu_relationship}
\sigma(z) = \mu^2(z,0), \quad z \in \Om.
\end{equation}

In the regularized method, the functions $\psi(z,k)$ in Step 1, $t(k)$ in Step 2, and $\mu(z,k)$ in Step 3 are computed for complex frequencies $|k|\leq R$, where $R$ serves as a regularization parameter \cite{Knudsen09} dependent upon the noise level.  In an ideal setting with no noise present, the reconstruction converges to the actual conductivity as $R\rightarrow\infty$ pointwise \cite{Knudsen09}.  
Our {\em a priori} method will take advantage of this fact, balancing the ability to choose $R$ large for ideal priors with fidelity to the data.


\section{Outline of the a priori method} \label{APriori_outline}
To motivate the \emph{a priori} scheme, we first note that the scattering transform can be written in terms of a scattering transform computed from a prior known conductivity distribution $\sigpr$ in the form $\T(k) = \Tpr+\mbox{perturbation}$.  To this end, given $\sigpr$, let $\Lampr$ denote the DN map corresponding to $\sigpr$, let $\psipr$ denote the CGO solution satisfying 
\begin{equation} \label{psipr_BIE}
\psipr(z,k)\vert_{\DOm} = e^{ikz}|_{\DOm}-\int_{\DOm}G_k(z-\zeta)(\Lampr-\Lambda_1)\psipr(\cdot,k) ds(\zeta),
\end{equation}
let $\Tpr$ denote the scattering transform satisfying
\begin{equation} \label{tpr_BIE}
\Tpr(k) = \int_{\DOm}e^{i\bar{k}\bar{z}}(\Lampr-\Lambda_1)\psipr(z,k) ds(z),
\end{equation}
and let $\mupr$ denote the CGO solution satisfying
\begin{equation}\label{Dbar_int_eqn_mupr}
\mupr(z,k) = 1 + \frac{1}{(2\pi)^2} \int_{\R^2} \frac{\Tpr(k')}{\bar{k'}(k-k')} e^{-i (zk' + \bar{z}\bar{k'})} \overline{\mupr(z,k')} \; dk'.
\end{equation}
These equations are valid for $\sigpr\in L^{\infty}(\Om)$ by \cite{Astala06}.  Subtracting \eqref{tpr_BIE} from \eqref{scatt_trans2}, we see
\begin{eqnarray*}
\T(k) -\Tpr(k) &=& \int_{\DOm}e^{i\bar{k}\bar{z}}(\Lambda_{\sigma}\psi -\Lambda_1\psi -\Lampr\psipr +\Lambda_1\psipr) ds\\
&=& \int_{\DOm}e^{i\bar{k}\bar{z}}(\Lambda_{\sigma}(\psi -\psipr) -\Lambda_1(\psi -\psipr) +(\Lambda_{\sigma} -\Lampr)\psipr) ds.
\end{eqnarray*}
Thus,
\begin{equation} \label{t_tpr}
\T(k) = \Tpr(k) + \int_{\DOm}e^{i\bar{k}\bar{z}}[(\Lambda_{\sigma}-\Lambda_1)(\psi -\psipr) +(\Lambda_{\sigma} -\Lampr)\psipr] ds.
\end{equation}
Formula \eqref{t_tpr} suggests the following scheme.  Given $\sigpr$, compute $\Lampr$ from a numerical forward solver, such as FEM, compute $\psipr$ and $\Tpr$ from \eqref{psipr_BIE} and \eqref{tpr_BIE}, respectively, compute $\T(k)$ from \eqref{t_tpr}, and use this $\T(k)$ in the D-bar method.  However, this natural approach has several drawbacks when applied to noisy data.  First of all, since the measured data $\Lambda_{\sigma}$ has noise, it is necessary to compute $\T(k)$ on a truncated domain $|k|\leq R$.  This means finer details encoded in large $|k|$ values of the prior will be lost.  Second, the numerical computation of $\Lampr$ itself introduces error that is not necessarily a good match to the noise in $\Lambda_{\sigma}$.  Thus, the term
$$\int_{\DOm}e^{i\bar{k}\bar{z}}(\Lambda_{\sigma} -\Lampr)\psipr ds$$ 
is not an accurate perturbation, and errors in $\Lambda_{\sigma}-\Lampr$ will be amplified by the exponentially growing functions $e^{i\bar{k}\bar{z}}$ and $\psipr$.

An alternative approach motivated by \eqref{t_tpr} is to define an approximation to the scattering transform piecewise by 
\begin{equation} \label{t_R1_R2}
\TRs(k) := \begin{cases} \T(k), &|k| \leq R_1 \\ \Tpr(k), &R_1 < |k| \leq R_2 \\ 0, &|k| > R_2 \end{cases}.
\end{equation}
In this approximation, the perturbation term in \eqref{t_tpr} is neglected for $|k| > R_1$, and the entire scattering transform is truncated for some $R_2 \geq R_1$.  Neglecting this term for $|k| > R_1$ is motivated by the fact that the size of $R_1$ is limited since $\T(k)$ will inevitably blow up for larger values of $|k|$  in the presence of noisy data. Since the computation of $\Tpr$ is noise-free and in general much more numerically robust than the computation of $\T$, we may select $R_2$ to be significantly larger than $R_1$. The larger the value of $R_2$, the stronger the influence of the \emph{a priori} information. 
The next question is then how to compute $\TRs$.  For $|k|\leq R_1$, $\T(k)$ can be computed using \eqref{scatt_trans2}, which is the same as computing $\TRBIE$, using the notation in \cite{Knudsen09}.  To avoid the problems that arise from computing $\Lampr$, the scattering transform $\Tpr$ can be computed directly from the definition of the scattering transform, provided $\sigpr\in C^2$.  Then, defining $\qpr = \frac{\Delta \sqrt{\sigpr}}{\sqrt{\sigpr}},$ the scattering transform $\Tpr$ is defined to be the nonlinear Fourier transform of $\qpr$ \cite{Nachman96}
\begin{equation} \label{tpr}
\Tpr(k) := \int_{\R^2} e^{i \bar k \bar z}\qpr(z) \psipr(z,k) dz, 
\end{equation}
where $\psipr$ is the solution of the Schr\"odinger equation \eqref{Schrod_eqn}, with $q = \qpr$. 
Once the scattering transform has been computed, the CGO solution $\mu$ can be solved from \eqref{Dbar_int_eqn}.  We define $\mutilde$ as the solution to
\begin{equation}\label{Dbar_int_eqn_tilde}
\mutilde(z,k) = 1 + \frac{1}{(2\pi)^2} \int_{|k|\leq R_2} \frac{\TRs(k')}{\bar{k'}(k-k')} e^{-i (zk' + \bar{z}\bar{k'})} \overline{\mutilde(z,k')} \; dk'.
\end{equation}
However, there is one more thing to note.  The Green's function for the D-bar operator $\dbar_k$ is $\frac{1}{\pi k}$, and so the solution \eqref{Dbar_int_eqn} to \eqref{Dbar_eq} is obtained from
\begin{equation}\label{Dbar_int_eqn_limit}
\mu(z,k) =  \lim_{R\rightarrow\infty}\left\{ \frac{1}{\pi R^2} \int_{|k|\leq R} \mu(z,k) dk + \frac{1}{(2\pi)^2} \int_{|k|\leq R} \frac{\T(k')}{\bar{k'}(k-k')} e^{-i (zk' + \bar{z}\bar{k'})} \overline{\mu(z,k')} \; dk' \right\},
\end{equation}
where the first term tends to $1$ as $R\rightarrow\infty$.  Thus,
\begin{equation}\label{Dbar_int_eqn_limit2}
\mutilde(z,k) \approx \frac{1}{\pi R_2^2} \int_{|k|\leq R_2} \mu(z,k) dk + \frac{1}{(2\pi)^2} \int_{|k|\leq R_2} \frac{\TRs(k')}{\bar{k'}(k-k')} e^{-i (zk' + \bar{z}\bar{k'})} \overline{\mutilde(z,k')} \; dk' .
\end{equation}
In practical computations, since $\mutilde(z,k)$ is unknown, the first integral is replaced by $1$, that is, its limit as $R_2\rightarrow\infty$.  Approximating $\mutilde$ in this term by $\mupr$, we have derived an equation for the approximation
\begin{equation}\label{Dbar_int_eqn_muR2}
\mu_{R_2}(z,k) = \frac{1}{\pi R_2^2} \int_{|k|\leq R_2} \mupr(z,k) dk + \frac{1}{(2\pi)^2} \int_{|k|\leq R_2} \frac{\TRs(k')}{\bar{k'}(k-k')} e^{-i (zk' + \bar{z}\bar{k'})} \overline{\mu_{R_2}(z,k')} \; dk'.
\end{equation}
If the prior coincides with the correct conductivity distribution, this method converges and introduces no artifacts as $R_1, R_2\rightarrow\infty$ by the convergence proof for the regularized D-bar method \cite{Knudsen09}.

The strength of the prior, or its influence on the reconstruction, depends on $R_2$.  If $R_2 = R_1$, then $\TRs(k) = \T_{\mbox{\tiny $R_1$}}^{\mbox{{\tiny \rm BIE}}}(k)$, and the only influence of the prior on the reconstruction is in the term
\begin{equation} \label{muint}
\muint(z) := \frac{1}{\pi R_2^2} \int_{|k|\leq R_2} \mupr(z,k) dk.
\end{equation}
We can exert control over the amount of influence this term has by introducing a weighting parameter $\alpha\in [0,1]$ and writing
\begin{equation}\label{Dbar_int_eqn_muR2_alpha}
\mualpha(z,k) =\alpha +(1-\alpha)\muint(z) + \frac{1}{(2\pi)^2} \int_{|k|\leq R_2} \frac{\TRs(k')}{\bar{k'}(k-k')} e^{-i (zk' + \bar{z}\bar{k'})} \overline{\mualpha(z,k')} \; dk',
\end{equation}
which is equivalent to \eqref{Dbar_int_eqn_tilde} if $\alpha=1$ and \eqref{Dbar_int_eqn_muR2} if $\alpha=0$.

\subsection{Computational considerations}

We now describe the numerical details for the computation of the conductivity distribution $\sigalpha$ corresponding to the CGO solution $\mualpha$ to (\ref{Dbar_int_eqn_muR2_alpha}), including numerical approximations for the necessary intermediate operators and functions. 

As described in \cite{MuellerBook}, we compute a finite-dimensional matrix approximation  ${\bf L}_\sigma$ to the DN map $\Lambda_\sigma$ by first computing the discrete ND map ${\bf R}_{\sigma}$ and then forming its inverse: ${\bf L}_\sigma = { \bf R}_\sigma^{-1}$. Denote  the number of linearly independent current patterns by $N$ and  the number of electrodes by $L$. To compute ${\bf R}_\sigma$, first orthonormalize the matrix of bipolar current patterns to obtain the set $\{J^m_l\}$, $m =1, \dots, N$ and $l = 1, \dots, L$, and then apply the appropriate change-of-basis formula to the voltages, yielding $\{V^m_l\}$. Then the ND matrix can be approximated by
$$
{\bf R}_\sigma(m,n) \approx  \sum_{l=1}^L \frac{\Delta \theta }{A }J^m_l V^n_l = \frac{\Delta \theta}{A} {\bf V}^\top {\bf J},
$$
where $A$ is the area of an electrode and $\Delta \theta$ is the angular distance between electrodes.

The discrete matrix approximation ${\bf L}_1$ of the DN map $\Lambda_1$ corresponding to homogeneous conductivity is obtained by first numerically solving the forward conductivity problem, using FEM for example, to create simulated voltage data for the case where $\sigma \equiv 1$ in $\Om$. The method described above for the computation of ${\bf L}_\sigma$ may then be used to compute ${\bf L}_1$ from this simulated data. 

The CGO solution $\psi|_{\DOm}$ is found by numerically solving the boundary integral equation \eqref{BIE} at the center $z_l$ of each electrode, as described in \cite{Herrera15}. As in \cite{Astala11}, we express the Faddeev Green's function $G_k$ as
$$
G_k(z) = \frac{1}{4 \pi} \Re( {\rm EI}(ikz)),
$$
where ${\rm EI}(z)$ is the exponential integral function, which in Matlab can be computed easily using the built-in function: ${\rm EI}(z) \approx 2* \mbox{\footnotesize{EXPINT}}(z)$. To form a matrix approximation $ {\bf \Gamma_k}$ for $G_k(z_l - \zeta_{l'})$, we must be careful of the logarthmic singularity that occurs when $l' = l$. We therefore discretize the surface of each electrode into $S$ points $z_{l_s}$, $s=1, \dots, S$, and compute
$$
{\bf \Gamma_k}(l,l') = \begin{cases}  \frac{A}{2 \pi} \Re(\mbox{\footnotesize{EXPINT}}(ik(z_{l}- \zeta_{l'}))), & l' \neq l \\ \frac{1}{2 \pi(S-1)} \sum_{s=1}^S \Re(\mbox{\footnotesize{EXPINT}}(ik(z_{l}- \zeta_{l_s}))) & l' = l \end{cases}.
$$

Denote by ${\bf b_k} = (b_1(k), \dots, b_N(k))^{\top}$ and ${\bf c_k}$ the vectors of coefficients for the functions $\psi(z,k)|_{\DOm}$ and $e^{ikz}|_{\DOm}$, respectively, expanded in the basis of orthonormalized current patterns, so that $\psi(z, k)|_{\DOm} \approx {\bf J} {\bf b_k}, e^{ik z}|_{\DOm}  \approx {\bf J} {\bf c_k}$. We may then approximate the convolution of $G_k$ with $(\Lambda_\sigma - \Lambda_1)\psi$  for each $z = z_l$ as a finite-dimensional vector:
$$
 \int_{\DOm} G_k(z - \zeta)(\Lam_\sigma - \Lam_1)\psi(\cdot,k) ds(\zeta) \approx {\bf \Gamma_k}{\bf J}( {\bf L}_\sigma - {\bf L}_1 ) {\bf b_k},
$$
and the discrete version of \eqref{BIE} is
$$
{\bf J b_k} = {\bf J c_k } - {\bf \Gamma_k J}( {\bf L}_\sigma- {\bf L}_1 ) {\bf b_k }. 
$$
Multipying through by the transpose of the orthonormal matrix ${\bf J}$ yields the linear system
\begin{equation} \label{BIE_lin_sys}
[ {\bf I} + {\bf J}^\top {\bf \Gamma_k J}( {\bf L}_\sigma- {\bf L}_1 ) ]{\bf b_k } = {\bf c_k},
\end{equation}
where ${\bf I}$ is the $N \times N$ identity matrix. In our implementation, this system was solved in Matlab using the {\footnotesize MLDIVIDE} function. 

The computation of the CGO solution $\psipr$ corresponding to the prior is handled quite differently. From \cite{Nachman96}, we know $\mupr$ satisfies the Lippmann-Schwinger type equation
\begin{equation} \label{LS_eqn_pr}
\mupr = 1 - g_k\ast (\qpr\mupr).
\end{equation}
where $g_k(z)=e^{-ikz} G_k(z)$ is a fundamental solution of the operator $-\Delta - 4ik\dbar$.
The numerical solution of \eqref{LS_eqn_pr} is based on ideas presented in \cite{Vainikko00}, and a complete description of the computational steps can be found in \cite{MuellerBook}. In short, we may write \eqref{LS_eqn_pr} as the linear system
\begin{equation} \label{LS_lin_sys}
[I + g_k \ast \qpr(\cdot)]\mupr = 1,
\end{equation}
which is solved for $\mupr$ using a matrix-free method such as {\footnotesize GMRES} \cite{Saad86} or {\footnotesize BICGSTAB} \cite{Vandervorst92}, which was used in our implementation, separating real and imaginary parts as required for the linear solver. The action of the linear operator on the left-hand side of \eqref{LS_lin_sys}  may be approximated efficiently using {\footnotesize FFT} and  {\footnotesize IFFT} operations. 

Once we obtain $\psipr$, the computation of $\Tpr$ from its definition \eqref{tpr} is accomplished using simple numerical quadrature over the mesh of $z$-values. The integral $\muint$ is likewise found by applying numerical quadrature to the integrand $\mupr$.

To obtain $\mualpha(z)$, we must solve the equation \eqref{Dbar_int_eqn_muR2_alpha} for each $z$, which involves modification of the computational methods described in \cite{MuellerBook}. We write \eqref{Dbar_int_eqn_muR2_alpha} as the linear system
\begin{equation} \label{Dbar_lin_sys}
[I - \mathcal A T(\bar \cdot)]\mualpha(z) = \alpha +(1-\alpha)\muint(z)
\end{equation}
where the actions of the operators $T$ and $\mathcal A$ are defined by
$$
Tf(k) = \frac{\TRs(k)}{4 \pi \bar k} e^{-i (zk' + \bar{z}\bar{k'})}  f(k), \quad  \quad \mathcal A g(k)  = \frac{1}{\pi} \int_{|k|\leq R_2}\frac{g(k)}{k-k'}dk'.
$$
Note that \eqref{Dbar_lin_sys} is not complex-linear due to the presence of the conjugate operator, so it is necessary to solve real and imaginary parts separately. This system is solved by again using a matrix-free solver such as {\footnotesize BICGSTAB}, where the action of $\mathcal A$ can be approximated by {\footnotesize FFT} and  {\footnotesize IFFT} operations. 

From $\mualpha$ we obtain the resulting conductivity distribution $\sigalpha = \mualpha^2(z,0)$, which is a result of both the EIT data $\Lam_\sigma$ and the \emph{a priori} data encoded into $\TRs$ and $\muint$.

\subsection{Constructing the prior: Blind Estimate Method}\label{BlindEstimate}
We now discuss the first of two possible methods for assigning conductivity values to the \emph{a priori} conductivity distribution $\sigpr$. Each of these two methods describes the construction of a discontinuous \emph{a priori} distribution $\tsigpr$. To satisfy the requirement that $\sigpr \in C^2(\Om)$, we later mollify $\tsigpr$ to obtain $\sigpr$. 

\begin{figure}[htb]
  \centering
    \includegraphics[width=0.45\textwidth]{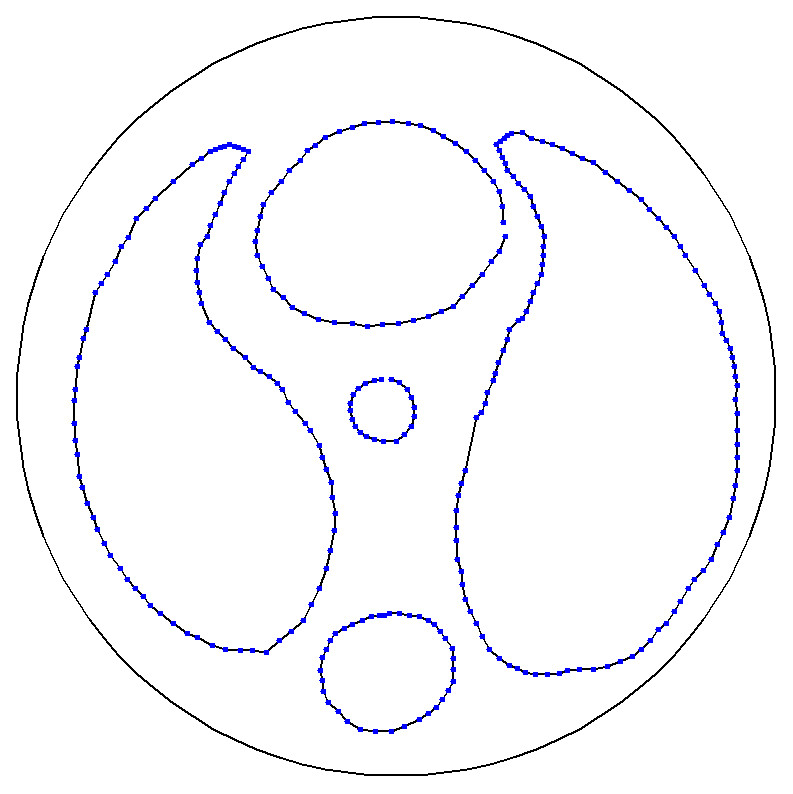}
  \caption{Simulated organ boundaries representing heart, lungs, aorta, and spine within a circular domain, used as the simulated \emph{a priori} information in our experiments. Organ boundaries are approximated by polygonal regions; the dots represent polygon vertices.}
\label{prior_boundaries}
\end{figure}

In both methods, knowledge is assumed of the spatial locations of boundaries for various domain inclusions (such as heart, lungs, etc.\ for thoracic imaging) in the plane of the electrodes. In a clinical setting, this could be obtained by extracting the organ boundaries from a CT scan to obtain polygonal approximations to the actual organ boundaries. In our simulations we created polygonal boundaries representing heart, lungs, aorta, and spine within a circular domain, as shown in Figure~\ref{prior_boundaries}. 

In the blind estimate method for assigning conductivity values to the prior, we simply make educated guesses for the conductivity values within each approximate organ boundary. These values can, for example, be estimated from literature sources wherein conductivity values for human tissue have been reported.  Let $P \subset \Om$ denote the polygonal region inside a particular approximate organ boundary and let $\{z_r\}$  be the finite set of points in the $z$-mesh used to construct $\sigpr$.  An approximate constant conductivity value $\sigma_{\mbox{\tiny \it P}}$ is selected for $P$, and we assign $\tsigpr(z_n) = \sigma_{\mbox{\tiny \it P}}$ for all $z_n \in \{z_r\} \cap P$. We repeat this process for all organ boundaries used in the prior, obtaining the conductivity distribution $\tsigpr$. Refinements to this process could be made by specifying regions of differing conductivities within individual organs if known inhomogeneities exist.

The blind estimate method is much computationally simpler and faster than the alternative extraction method which will be described next. However, if any pathologies have developed between the time of the initial CT scan and the time of the EIT scan, these pathologies will not be reflected in the prior, and their expression in the final reconstruction is therefore entirely dependent on the EIT data. The full \emph{a priori} scheme with the blind estimate method for constructing the prior is outlined in Algorithm~\ref{BlindEstimateAlg}. 


\begin{algorithm}[h]
\floatname{algorithm}{{\footnotesize ALGORITHM}}
\begin{algorithmic}[1]
{\footnotesize
\Proc{Obtain \emph{a priori} information and EIT data:}
	\DashState Form polygonal approximations to organ boundaries.
	\DashState Collect EIT data and compute ${\bf L}_\sigma$.
	\DashState Use FEM to simulate homogeneous data and compute ${ \bf L}_1$. 
\EndProc
\State Form computational grids for the $k$ and $z$ planes. 
\State Make blind estimates for conductivity values to form $\sigpr(z)$. 
\Proc{Compute $\sigalpha$:}
	\DashState Select $R_1, R_2$ with $R_1 \leq R_2$. 
	\DashState Compute $\qpr = \Delta \sqrt{\sigma} / \sqrt{\sigma}$.  
	\For{$|k| \leq R_1$}
		\DashState Solve \eqref{BIE_lin_sys} for ${\bf b_k}$ to get $\psi|_{\DOm} \approx {\bf J b_k}$
		\DashState Compute $\T(k)$ from \eqref{scatt_trans2}. 
	\EndFor
	\For{$R_1 < |k| \leq R_2$}
		\DashState Solve \eqref{LS_lin_sys} for  $\mupr = e^{-ikz}\psipr$. 
		\DashState Compute $\Tpr(k)$ from \eqref{tpr}.	
	\EndFor
	\For{$z \in \Om$}
		\DashState Solve \eqref{Dbar_lin_sys} for $\mualpha(z,\cdot)$.
		\DashState Compute $\sigalpha(z)$ from \eqref{sigma_mu_relationship}. 
	\EndFor 
\EndProc
\caption{\emph{A priori} scheme, Blind Estimate Method} \label{BlindEstimateAlg}
}
\end{algorithmic}
\end{algorithm}

\subsection{Constructing the prior: Extraction Method}\label{Extraction}

In the extraction method, we first compute a reconstruction $\sigma$ from the EIT data alone using the D-bar method described in \S\ref{Background}.  Note that the first-order approximation $\Texp$ to the scattering transform as given in \cite{Isaacson04} could be used here as well to obtain an initial reconstruction.  We then extract conductivity values from this reconstruction to obtain estimated values for the prior. This method has advantages over the blind estimate method in that pathologies not present in the CT scan data but that are apparent in the reconstruction $\sigma$ may be included in the prior. The full \emph{a priori} scheme with the extraction method for constructing the prior is outlined in Algorithm~\ref{ExtractionAlg}. 

In our experiments with simulated data, we developed and used the following techniques for extracting approximate conductivity values from the reconstruction $\sigma$ to be assigned to the lungs, heart, aorta, spine, and background in the prior. In what follows, we denote by $\{z'_s\}$ the finite set of points in the $z$-mesh used to construct $\sigma$, and the polygonal regions inside specific organ boundaries by $P_{\mbox{\tiny heart}}$, $P_{\mbox{\tiny spine}}$, etc. 

\begin{enumerate}
\item {\bf Lungs.} We examine each lung in the $\sigma$ reconstruction and compare the appearance of the lungs to the \emph{a priori} approximate lung boundaries. If, based on the  reconstruction $\sigma$, the lungs appear to be free of pathology (i.e.\ there are no suspicious inhomogeneities within the regions $P_{\mbox{\tiny l lung}}$ and $P_{\mbox{\tiny r lung}}$), then the following method may be used. Find the set $\tilde{P}_{\mbox{\tiny l lung}} := \{z'_s\} \cap P_{\mbox{\tiny l lung}}$, and compute the average value
\begin{equation}\label{lung_comp}
\sigma_{\mbox{\tiny l lung}} := \frac{1}{M} \sum_{m = 1}^M \sigma(z'_m), \quad z'_m \in \tilde P_{\mbox{\tiny l lung}},
\end{equation}
where $M$ denotes the number of points $z'_m$ in $\tilde{P}_{\mbox{\tiny l lung}}$, and then assign $\sigpr(z_n) = \sigma_{\mbox{\tiny l lung}}$ for all $z_n \in \{z_r\} \cap {P}_{\mbox{\tiny l lung}}$. This process is then repeated for the right lung. On the other hand, if the  reconstruction $\sigma$ reveals possible lung pathologies in the form of inhomogeneities within a lung region, then the method can be revised in the following way. Assuming (without loss of generality) that one or more inhomogeneities appear in the left lung, divide the region $P_{\mbox{\tiny l lung}}$ into a finite number of connected subsets $S_j \subset P_{\mbox{\tiny l lung}}$, where each subset represents an area of fairly homogeneous conductivity in the  reconstruction $\sigma$. Then for each $S_j$, compute the average conductivity over the points $z'_m \in \{z'_s\} \cap S_j$ and assign this value to $\tsigpr(z_n)$ for $z_n \in \{z_r\} \cap S_j$. 

\item {\bf Heart and aorta.} To compute $\tsigpr$ values for the heart region, one could potentially employ a similar method to that described for the lung regions. However, the position and shape of the reconstructed heart is more sensitive to noise level and truncation radius than the lung, and therefore using the anatomical position in the prior may include extraneous pixels.  See Figure \ref{tBIE_recons} as an example of how the size and shape of the reconstructed heart can vary.  The aorta, on the other hand, is typically invisible in the  reconstruction $\sigma$, so such a method could not be used to assign $\tsigpr$ values within the aorta. The following method is therefore given as an alternative to the method used for the lungs.  First, define the quantities
$$
\sigma_{\mbox{\tiny max}} := \max_{z'_m \in \{z'_s\}\cap \Om} \{\sigma(z'_m)\}, \quad \quad \sigma_{\mbox{\tiny  min}} := \min_{z'_m \in \{z'_s\} \cap \Om} \{\sigma(z'_m)\},
$$
and compute the value $\tau = \sigma_{\mbox{\tiny  min}} + c (\sigma_{\mbox{\tiny max}} - \sigma_{\mbox{\tiny  min}})$ where $c \in (0.5, 1)$ is selected empirically.  A good choice for $c$ should optimally result in the set $H := \{ z'_m \in \{z'_s\}\cap \Om : \sigma(z'_m) \geq \tau \}$  being selected so as to be a connected subset of $\Omega$ and roughly the same size as the region $P_{\mbox{\tiny heart}}$, and such a $c$ may vary depending on noise levels and choice of truncation radius in the computation of $\sigma$. We find the set $H$ and compute
$$
\sigma_{\mbox{\tiny heart}} := \frac{1}{M} \sum_{m = 1}^M \sigma(z'_m), \quad z'_m \in H,
$$
where $M$ denotes the number of points in $H$. Finally, assign $\tsigpr(z_n) = \sigma_{\mbox{\tiny heart}}$  for all $z_n \in \{z_r\} \cap P_{\mbox{\tiny heart}}$, and, since the aorta likely has conductivity values very similar to those of the heart, further assign these same values to $\tsigpr$ in the region $P_{\mbox{\tiny aorta}}$. 

\item {\bf Spine.} Due to its small size, the reconstruction of the spine typically has very poor spatial resolution and its appearance and associated conductivity values can vary widely in the  reconstruction $\sigma$ in the presence of noise. Since we can usually assume that the spine is one of the most resistive objects in a thoracic EIT scan, we simply assign $\sigpr(z_n) = \sigma_{\mbox{\tiny min}}$ for all $z_n \in \{z_r\} \cap P_{\mbox{\tiny spine}}$. 

\item {\bf Background.} We define the background to be the set $P_{\mbox{\tiny bg}} := \Om -  \cup_j P_j$ where each $P_j$ corresponds to an organ boundary included in the prior, and assign values to $P_{\mbox{\tiny bg}}$ according to the following method. Compute the quantities $\tau_1 = \sigma_{\mbox{\tiny  min}} + c _1(\sigma_{\mbox{\tiny max}} - \sigma_{\mbox{\tiny  min}})$ and $\tau_2 = \sigma_{\mbox{\tiny  min}} + c_2 (\sigma_{\mbox{\tiny max}} - \sigma_{\mbox{\tiny  min}})$ where $c_1, c_2 \in (0, 1)$, $c_1 < c_2$. Find the set $B:= \{ z'_m \in \{z'_s\} : \sigma(z'_m) \in [\tau_1,\tau_2] \}$, and compute
$$
\sigma_{\mbox\tiny bg} := \frac{1}{M} \sum_{m = 1}^M \sigma(z'_m), \quad z'_m \in B,
$$
where $M$ denotes the number of points in $B$. Assign $\sigpr(z_n) = \sigma_{\mbox{\tiny bg}}$ for all $z_r \in \{z_n\} \cap P_{\mbox{\tiny bg}}$. As with the value $c$ used for the heart, the values $c_1$ and $c_2$ must be selected empirically, and may once again vary depending on noise levels and choice of truncation radius in the computation of $\sigma$. Since the lungs, which have low conductivity compared to the background, tend to dominate the reconstruction, it is usually advantageous to choose $c_1, c_2$ to be skewed to the upper end of the scale $(0,1)$. Well-chosen $c_1$ and $c_2$ should result in the set $B$ excluding most of the region corresponding to the lungs and spine, as well as the high conductivity region inside the heart. 
\end{enumerate}

\begin{algorithm}[h]
\floatname{algorithm}{{\footnotesize ALGORITHM}}
\begin{algorithmic}[1]
{\footnotesize
\Proc{Obtain \emph{a priori} spatial information and EIT data:}
	\DashState Form polygonal approximations to organ boundaries.
	\DashState Collect EIT data and compute ${\bf L}_\sigma$.
	\DashState Use FEM to simulate homogeneous data and compute ${ \bf L}_1$. 
\EndProc
\Proc{Compute conductivity $\sigma(z)$ using standard D-bar methods:} \label{begin_sigma_comp}
	\DashState Form computational grids for the $k$ and $z$ planes for both $\sigma$ and $\sigpr$. 
	\DashState Select a truncation radius $R_1$. 
	\For{$|k| < R_1$}
		\DashState  Solve \eqref{BIE_lin_sys} for ${\bf b_k}$ to get $\psi|_{\DOm} \approx {\bf J b_k}$
		\DashState Compute $\T(k)$ from \eqref{scatt_trans2}. 
	\EndFor
	\For{$z \in \Om$}
		\DashState Solve the ($R_1$-truncated) equation \eqref{Dbar_int_eqn} for $\mu(z,\cdot)$.
		\DashState Compute $\sigma(z)$ from \eqref{sigma_mu_relationship}. 
	\EndFor \label{end_sigma_comp}
\EndProc
\State Extract conductivity values from $\sigma(z)$ to form $\sigpr$.
\Proc{Compute $\sigalpha(z)$:} 
	\DashState Compute $\qpr = \Delta \sqrt{\sigma} / \sqrt{\sigma}$.
	\DashState Select $R_2  \geq R_1$.
	\For{$R_1 < |k| \leq R_2$}
		\DashState Solve \eqref{LS_eqn_pr} for $\psipr(\cdot, k)$. 
		\DashState Compute $\Tpr(k)$ from \eqref{tpr}.	
	\EndFor
	\DashState Form $\TRs$. 
	\DashState Select  $\alpha$.
	\For{$z \in \Om$}
		\DashState Solve \eqref{Dbar_lin_sys} for $\mualpha(z,\cdot)$.
		\DashState Compute $\sigalpha(z)$ from \eqref{sigma_mu_relationship}. 
	\EndFor 
\EndProc 
 \caption{\emph{A priori} scheme with Extraction Method} \label{ExtractionAlg}
}
\end{algorithmic}
\end{algorithm}

\subsection{Iterative approaches}\label{sec:Iteration}

The \emph{a priori} schemes described in the previous pages may be used alone to obtain a reconstruction $\sigalpha$, but there is potential for further refinement of spatial resolution through the use of iterative approaches. The motivation is to take advantage of the enhanced spatial resolution in $\sigalpha$ to construct a new prior that is more accurate than the original in terms of conductivity values and possible pathologies. We may include in this updated prior any new information that appears in the reconstruction $\sigalpha$. This may be especially advantageous if the blind estimate method was used to construct the original prior, but the patient has since developed some pathology that is visible in the EIT data. Another situation where iteration may provide enhanced results is if we desire to use the extraction method, but the reconstruction $\sigma$ has very poor spatial resolution. In \S\ref{Results} we provide an example of the first of these scenarios, using simulated data. 

The computational steps in these iterative approaches are the following: (1) obtain the reconstruction $\sigalpha$, using either the blind estimate or extraction method to assign conductivity values to the prior, (2) use the extraction method described in \S\ref{Extraction} to extract conductivity values from $\sigalpha$ (rather than from $\sigma$), (3)  use these extracted conductivity values to form an updated prior $\sigpr'$,  (4) repeat the \emph{a priori} scheme using the original EIT data with the updated prior $\sigpr'$ to obtain an updated reconstruction $\sigalpha'$. This entire process could potentially be repeated again if desired to obtain a second iterate $\sigalpha''$.

\section{Results} \label{Results}

We now present the results of our test problem using simulated data. In this test problem, we simulate a situation in which \emph{a priori} information is available from a previous CT scan, but the patient has since developed a pleural effusion in one lung. We tested the previously described \emph{a priori} scheme using both the blind estimate and extraction methods for assigning conductivity values, and an iteration step as described above was performed on the results from the reconstructions using the blind estimate method. 

\begin{figure}[ht]
  \centering
    \includegraphics[width=0.4\textwidth]{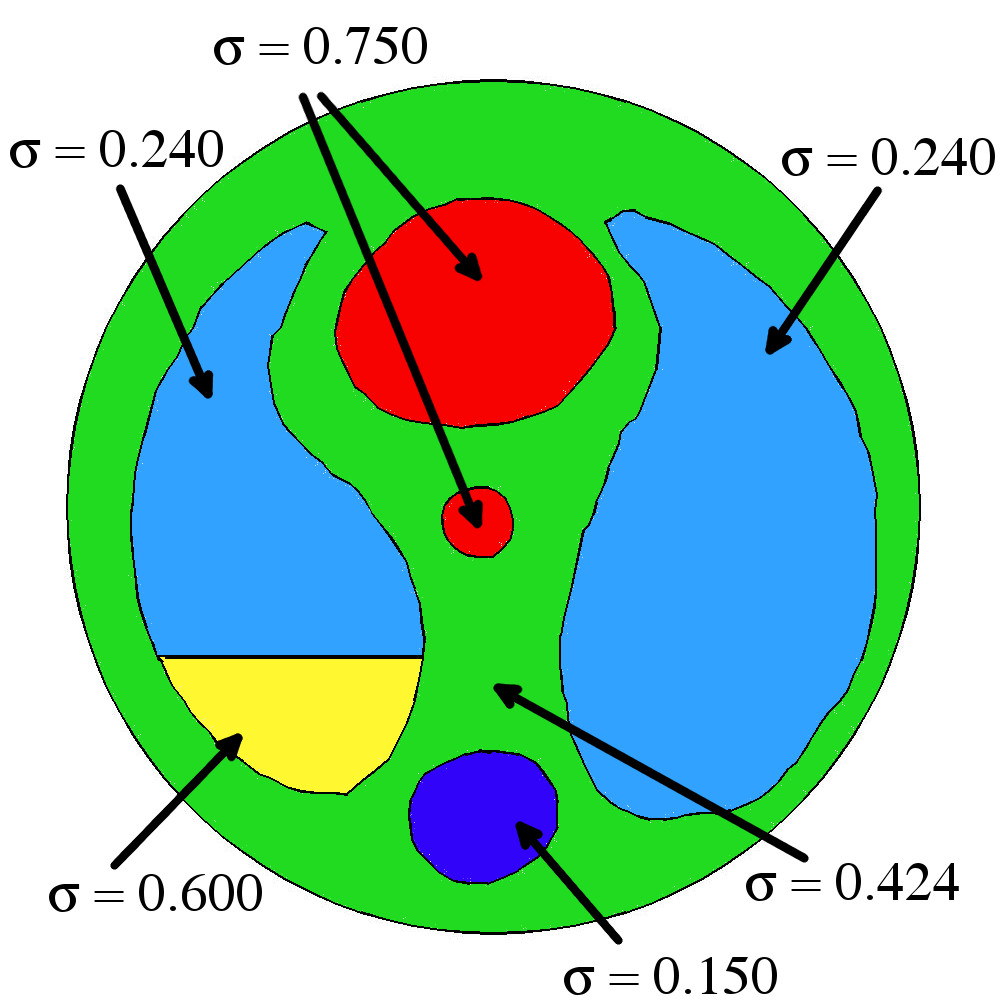}
  \caption{Phantom representing pleural effusion. Conductivity values are in S/m.}
\label{phantom_fluid_in_lung}
\end{figure}

\begin{figure}[ht]
  \centering
\begin{subfigure}{0.325\textwidth}
    \includegraphics[width=\textwidth]{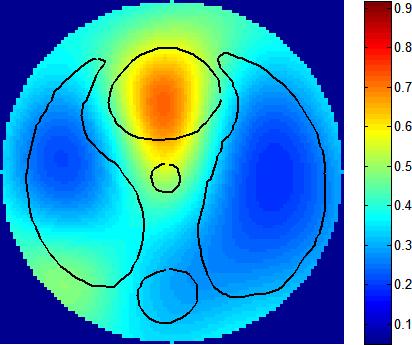}
  \caption{Noise level = 0\%.}
\end{subfigure}
\begin{subfigure}{0.325\textwidth}
    \includegraphics[width=\textwidth]{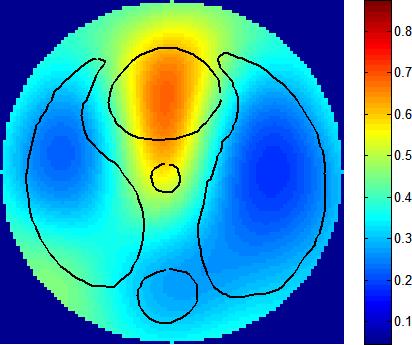}
  \caption{Noise level = 0.1\%.}
\end{subfigure}
\begin{subfigure}{0.325\textwidth}
    \includegraphics[width=\textwidth]{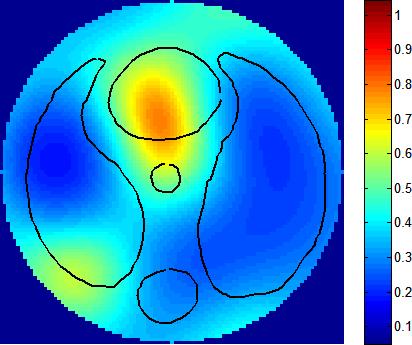}
  \caption{Noise level = 0.2\%.}
\end{subfigure}
\caption{Plots of the reconstructions $\sigma$ computed using the regularized D-bar method of \S\ref{Background} (see also lines \ref{begin_sigma_comp}--\ref{end_sigma_comp} of Algorithm~\ref{ExtractionAlg}) with truncation radius $R_1 = 3.8$, at noise levels 0\%, 0.1\%, and 0.2\%, with superimposed actual organ boundaries. Each noise case is plotted on its own scale; these scalings will be used for all plots concerning each noise level within this paper. }
\label{tBIE_recons}
\end{figure}

We assumed that the organ boundaries without pleural effusion were given by the polygonal approximations shown in Figure~\ref{prior_boundaries}, and we assumed a circular domain of radius 143.2 mm. We therefore created a phantom with these same organ boundaries, domain shape, and dimensions. Conductivity values were assigned to the phantom heart, lungs, aorta, and spine, and the FEM method with the complete electrode model including contact impedance described in \cite{Murphy09} was used to generate voltage data. To simulate a pleural effusion, conductivity was increased in the phantom in the bottom of the left lung. We will use the convention that the left lung appears on the left-hand side of the image. The phantom with assigned conductivity values is shown in Figure~\ref{phantom_fluid_in_lung}.  Random zero-mean Gaussian noise was added to the simulated voltages at 0\%, 0.1\%, and 0.2\% of the maximum voltage values; the D-bar reconstructions $\sigma$ using the method of \S\ref{Background} for each of the three noise cases are shown in Figure~\ref{tBIE_recons}.  All reconstructions in this section, including the reconstructions $\sigma$ and the results of the \emph{a priori} schemes, were computed using a $z$-mesh with $101 \times 101$ elements, $R_1 = 3.8$, and we tested values $R_2 \in \{3.8, 5.0, 7.5, 10\}$, and $\alpha \in \{ 0, 0.5, 0.75, 0.9\}$.  Note that the {\em a priori} organ boundaries are the correct boundaries, and the conductivity values and distributions will be modified.

\subsection{Blind estimate method applied to test problem} \label{BlindEstimateApplied}
\begin{figure}[ht]
  \centering
    \includegraphics[height = 0.88\textheight]{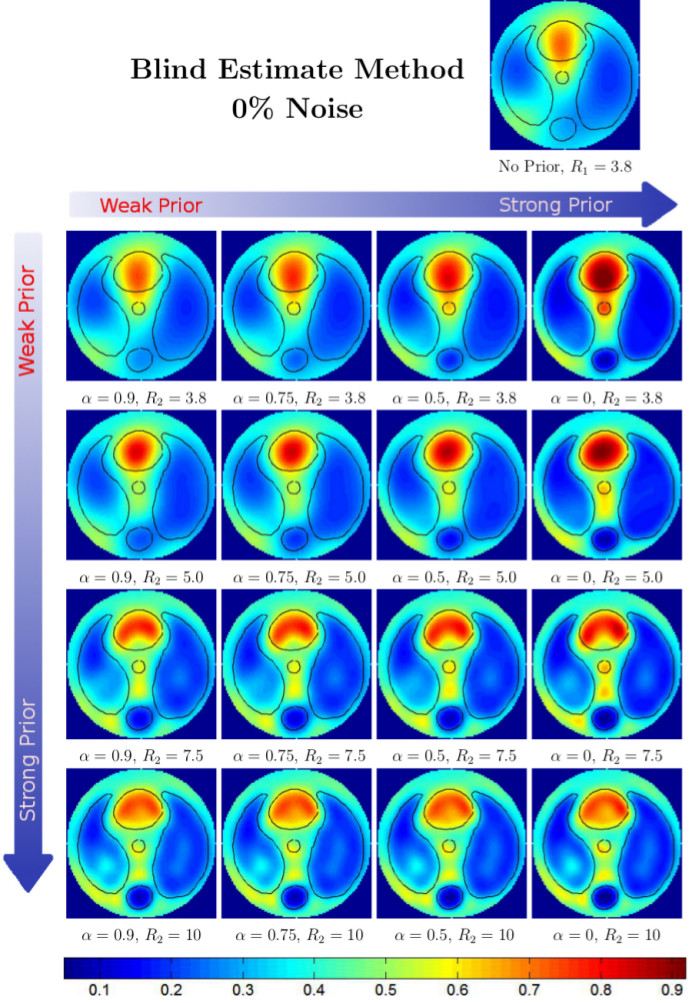}
  \caption{Reconstructions $\sigalpha$ for the 0\% noise case using the \emph{a priori} scheme with the blind estimate method (before iteration step), with various values of $\alpha$ and $R_2$. The reconstruction with no prior is at the top for comparison. The strength of the prior increases moving left to right and top to bottom. The scale bar at the bottom applies to all reconstructions at this noise level.}
\label{fig:BlindEst_NoNoise}
\end{figure}

\begin{figure}[ht]
  \centering
    \includegraphics[height = 0.88\textheight]{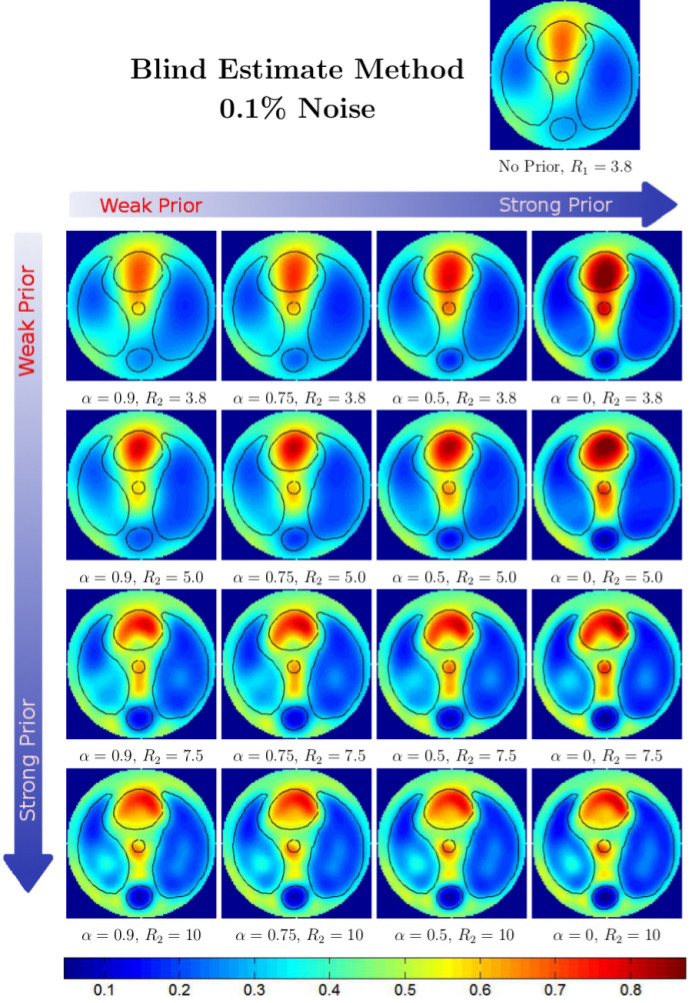}
  \caption{Reconstructions $\sigalpha$ for the 0.1\% noise case using the \emph{a priori} scheme with the blind estimate method (before iteration step), with various values of $\alpha$ and $R_2$. The reconstruction with no prior is at the top for comparison. The strength of the prior increases moving left to right and top to bottom. The scale bar at the bottom applies to all reconstructions at this noise level.}
\label{fig:BlindEst_0p1}
\end{figure}

\begin{figure}[ht]
  \centering
    \includegraphics[height = 0.88\textheight]{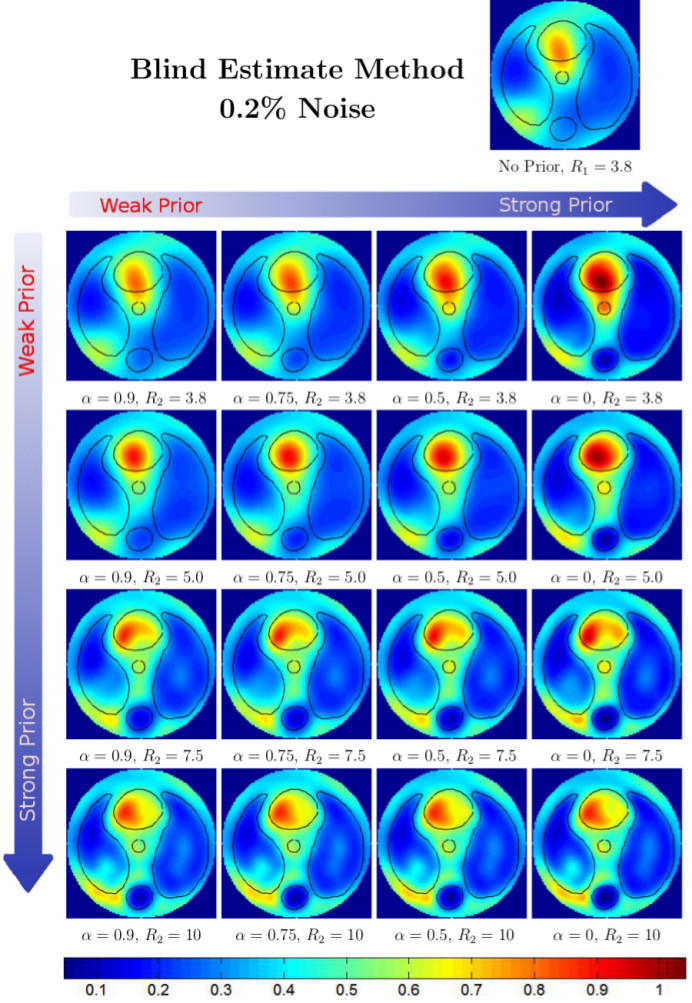}
  \caption{Reconstructions $\sigalpha$ for the 0.2\% noise case using the \emph{a priori} scheme with the blind estimate method (before iteration step), with various values of $\alpha$ and $R_2$. The reconstruction with no prior is at the top for comparison. The strength of the prior increases moving left to right and top to bottom. The scale bar at the bottom applies to all reconstructions at this noise level.}
\label{fig:BlindEst_0p2}
\end{figure}

We assigned ``blind estimate'' \emph{a priori} conductivity values representing a phantom with two homogeneous lungs with conductivity $0.200$ S/m, in contrast to the actual values displayed in Figure \ref{phantom_fluid_in_lung}.  The values for the background, heart, aorta, and spine differed slightly from the actual values.  These ``blind estimates'' were used for all three noise cases, and are given in Table~\ref{table:conductivities2}. Reconstructions using the blind estimate method can be seen in Figures~\ref{fig:BlindEst_NoNoise}, \ref{fig:BlindEst_0p1}, \ref{fig:BlindEst_0p2}.

Given the obvious lung pathology apparent in the reconstructions, we then performed an iteration step wherein the left lung was divided into two regions, which we shall refer to as the ``lung top'' and ``lung bottom,'' separated by a horizontal line segment.  We computed conductivity values for the iterate $\sigpr'$ separately in each of these two regions, using the methods described in \S\ref{Extraction} to extract conductivity values from the $\sigalpha$ reconstruction with $R_2 = 5.0$ and $\alpha = 0.75$. In the computation of the values for the heart, aorta, and background, we selected values $c = 0.85$, $c_1 = 0.25$, and $c_2 =0.95$; for simplicity, we used these same values in all noise cases. The resulting conductivity values used in $\sigpr'$ are also given in Table~\ref{table:conductivities2}. 

\begin{figure}[htb]
  \centering
    \includegraphics[width=0.5\textwidth]{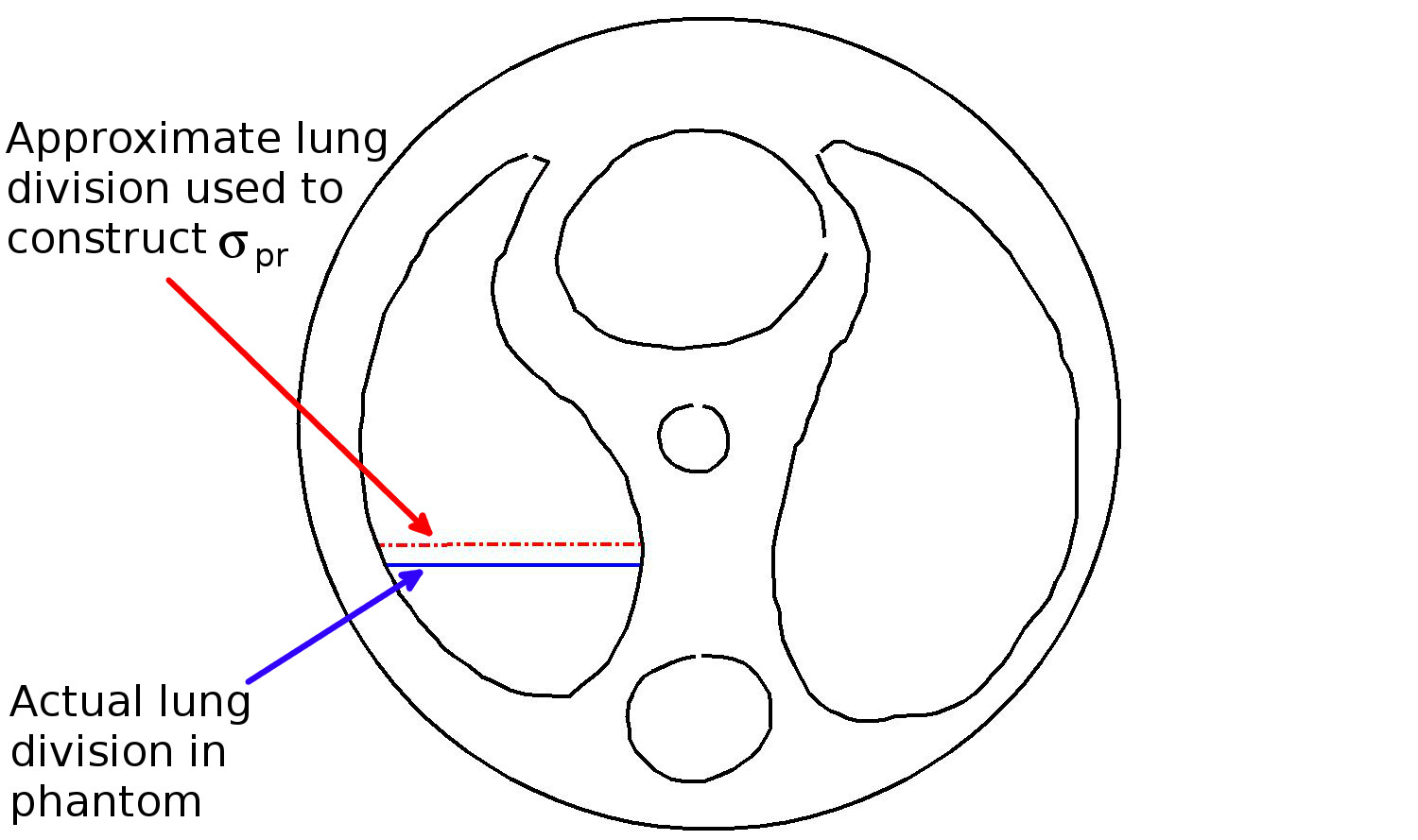}
  \caption{Locations of dividing line between the ``lung top'' and ``lung bottom.'' The dividing line used in the phantom is indicated by a solid line. The approximate dividing line (dashed line) was used in the extraction method and the iteration step for the blind estimate method, and was obtained by visually inspecting the $\sigalpha$ reconstructions from the blind estimate method.}
\label{lung_division}
\end{figure}

\begin{figure}[ht]
  \centering
    \includegraphics[height = 0.88\textheight]{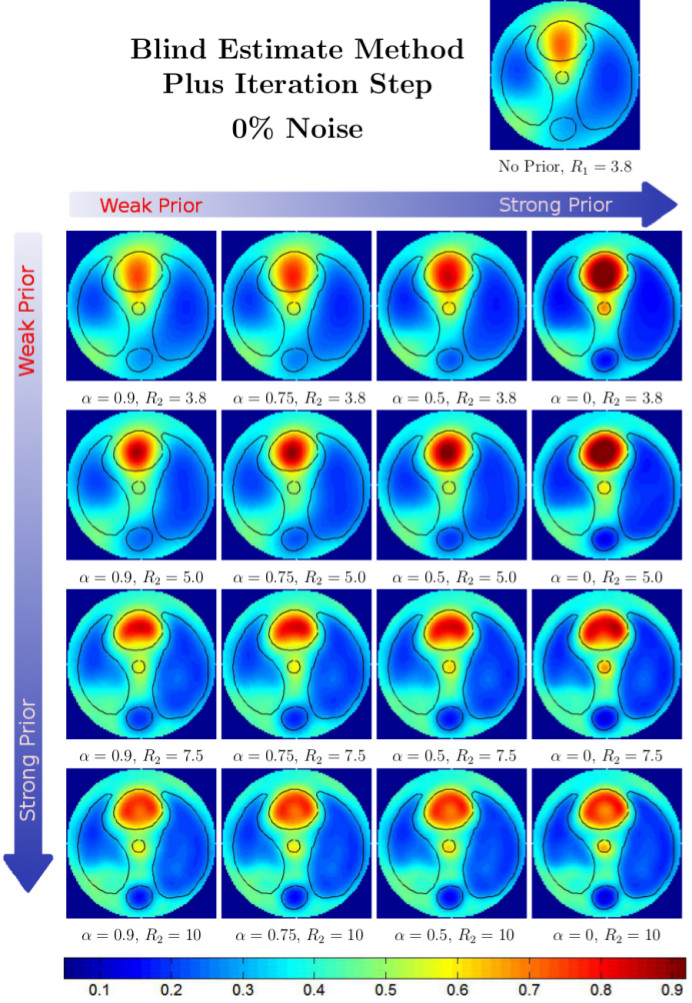}
  \caption{Reconstructions $\sigalpha'$ for the 0\% noise case using the \emph{a priori} scheme with the blind estimate method plus one iteration step, with various values of $\alpha$ and $R_2$. The reconstruction with no prior is at the top for comparison. The strength of the prior increases moving left to right and top to bottom. The scale bar at the bottom applies to all reconstructions at this noise level.}
\label{fig:BlindEst_Iter_NoNoise}
\end{figure}

\begin{figure}[ht]
  \centering
    \includegraphics[height = 0.88\textheight]{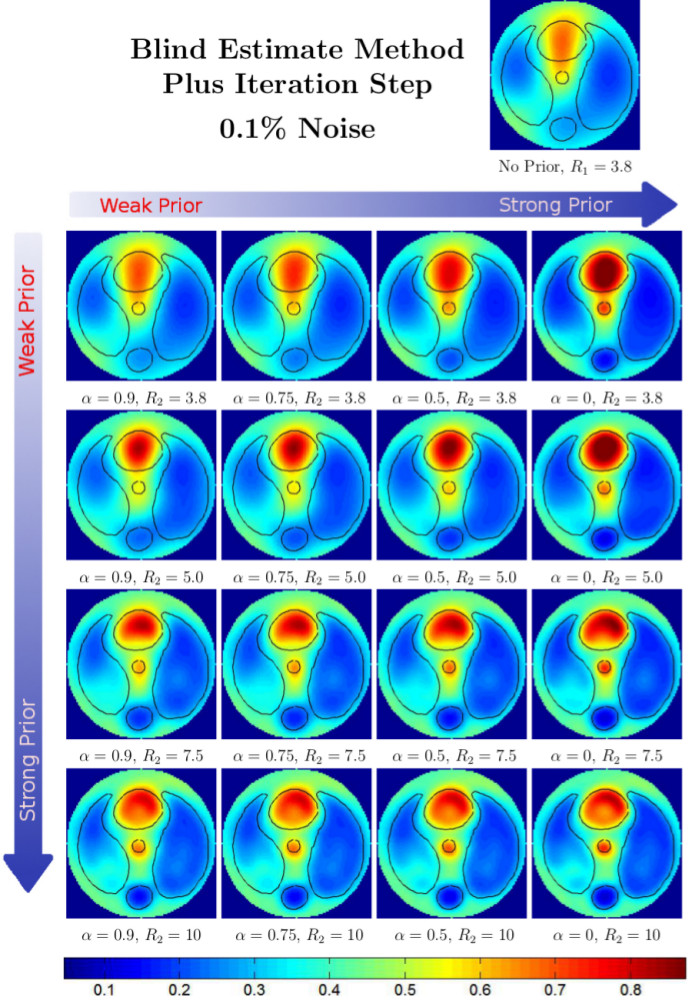}
  \caption{Reconstructions $\sigalpha'$ for the 0.1\% noise case using the \emph{a priori} scheme with the blind estimate method plus one iteration step, with various values of $\alpha$ and $R_2$. The reconstruction with no prior is at the top for comparison. The strength of the prior increases moving left to right and top to bottom. The scale bar at the  bottom applies to all reconstructions at this noise level.}
\label{fig:BlindEst_Iter_0p1}
\end{figure}

\begin{figure}[ht]
  \centering
    \includegraphics[height = 0.88\textheight]{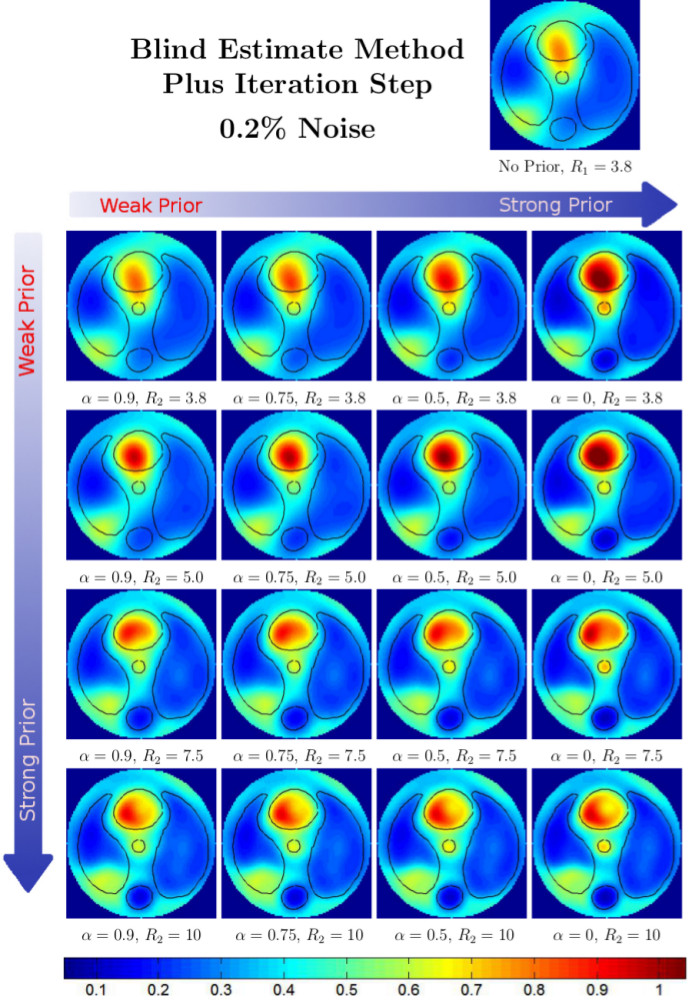}
  \caption{Reconstructions $\sigalpha'$ for the 0.2\% noise case using the \emph{a priori} scheme with the blind estimate method plus one iteration step, with various values of $\alpha$ and $R_2$. The reconstruction with no prior is at the top for comparison. The strength of the prior increases moving left to right and top to bottom. The scale bar at the bottom applies to all reconstructions at this noise level.}
\label{fig:BlindEst_Iter_0p2}
\end{figure}

The location of the dividing line between lung top and lung bottom was chosen by visually inspecting the $\sigalpha$ reconstructions and selecting a horizontal line at which to form the division.  For simplicity, we used the same approximate dividing line in all three noise cases. This approximate dividing line is compared to the actual lung division used to create the phantom in Figure~\ref{lung_division}.  The reconstructions resulting from the iteration step can be seen in Figures~\ref{fig:BlindEst_Iter_NoNoise}, \ref{fig:BlindEst_Iter_0p1}, and \ref{fig:BlindEst_Iter_0p2}.

\begin{figure}[htb]
  \centering
\begin{subfigure}{0.24\textwidth}
    \includegraphics[width=\textwidth]{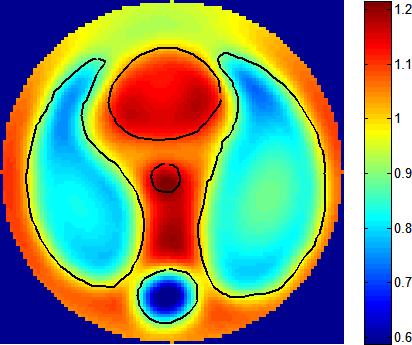}
  \caption{$R_2 = 3.8$}
\end{subfigure}
\begin{subfigure}{0.24\textwidth}
    \includegraphics[width=\textwidth]{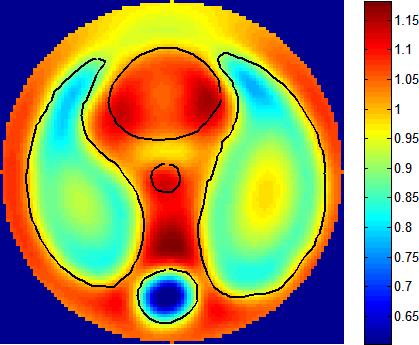}
  \caption{$R_2 = 5.0$}
\end{subfigure}
\begin{subfigure}{0.24\textwidth}
    \includegraphics[width=\textwidth]{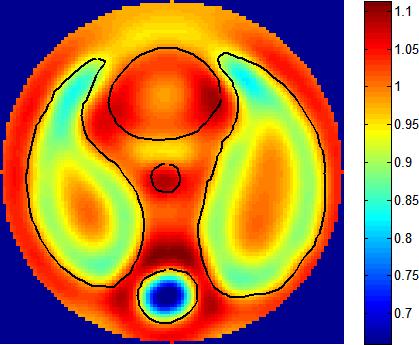}
  \caption{$R_2 = 7.5$}
\end{subfigure}
\begin{subfigure}{0.24\textwidth}
    \includegraphics[width=\textwidth]{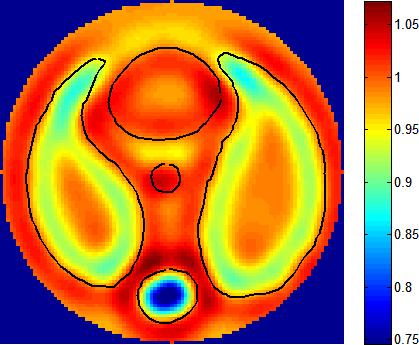}
  \caption{$R_2 = 10$}
\end{subfigure}
\caption{Plots of the real part of $\muint$ used in the simulations with various truncation radii $R_2$. Since $\muint \to 1$ as $R_2 \to \infty$, the scale must be adjusted for each value of $R_2$ for best viewing results.}
\label{muint_recons}
\end{figure}

It is evident from the reconstructions in Figures~\ref{fig:BlindEst_NoNoise}, \ref{fig:BlindEst_0p1}, \ref{fig:BlindEst_0p2} that with or without noise the blind estimate method without iteration detects the pleural effusion provided a very small value of $\alpha$ is not combined with a small value of $R_2$.  The case $\alpha=0$ and $R_2 = 3.8$ corresponds to using only $\muint$ as in equation \eqref{Dbar_int_eqn_muR2} and the scattering transform from the regularized D-bar method without a prior. Increasing $\alpha$ and decreasing $R_2$ weakens the influence of the prior.  In the case of a strong prior, the prior dominates the reconstruction in the blind estimate method, resulting in good spatial resolution of organ boundaries, but poor detection of the effusion.  Adding the iteration step described in \S\ref{sec:Iteration} results in excellent detection of the effusion in every case, and the aorta and spine can be clearly seen with excellent spatial accuracy even with a prior of medium weight, such as $\alpha = 0.5, R_2=7.5$.  The ringing effect seen in the reconstructed spine and lungs for small $\alpha$ or large $R_2$ is likely due to the influence of the term $\muint$ on the reconstruction.  This term provides excellent spatial resolution of the prior, detecting edges extremely well, but introducing ringing since it becomes increasingly uniform, tending to $1$ as $R_2\rightarrow\infty$.  This effect can be seen in Figure \ref{muint_recons}.

\subsection{Extraction method applied to test problem}
In assigning approximate conductivity values, for each noise level we first reconstructed $\sigma$ with $R_1 = 3.8$ using the D-bar method with no {\em a priori} information. The  reconstruction of $\sigma$ is plotted along with the \emph{a priori} organ boundaries in Figure~\ref{tBIE_recons}.

\begin{figure}[ht]
  \centering
    \includegraphics[height = 0.88\textheight]{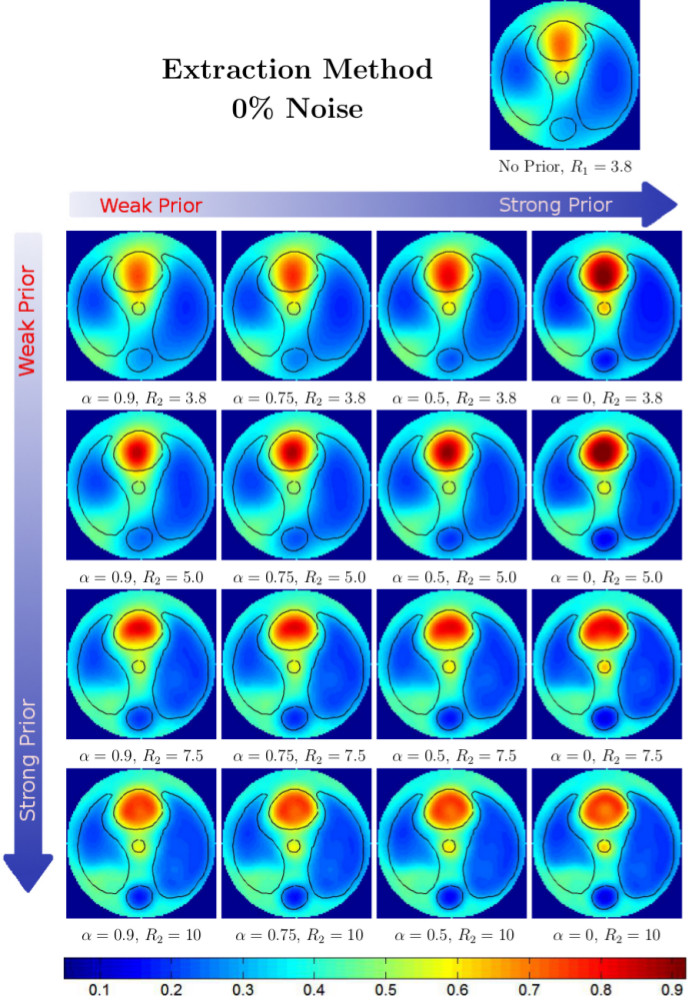}
  \caption{Reconstructions $\sigalpha$ for the 0\% noise case using the \emph{a priori} scheme with the extraction method, with various values of $\alpha$ and $R_2$. The reconstruction with no prior is at the top for comparison. The strength of the prior increases moving left to right and top to bottom. The scale bar at the bottom applies to all reconstructions at this noise level.}
\label{fig:Extraction_NoNoise}
\end{figure}

\begin{figure}[ht]
  \centering
    \includegraphics[height = 0.88\textheight]{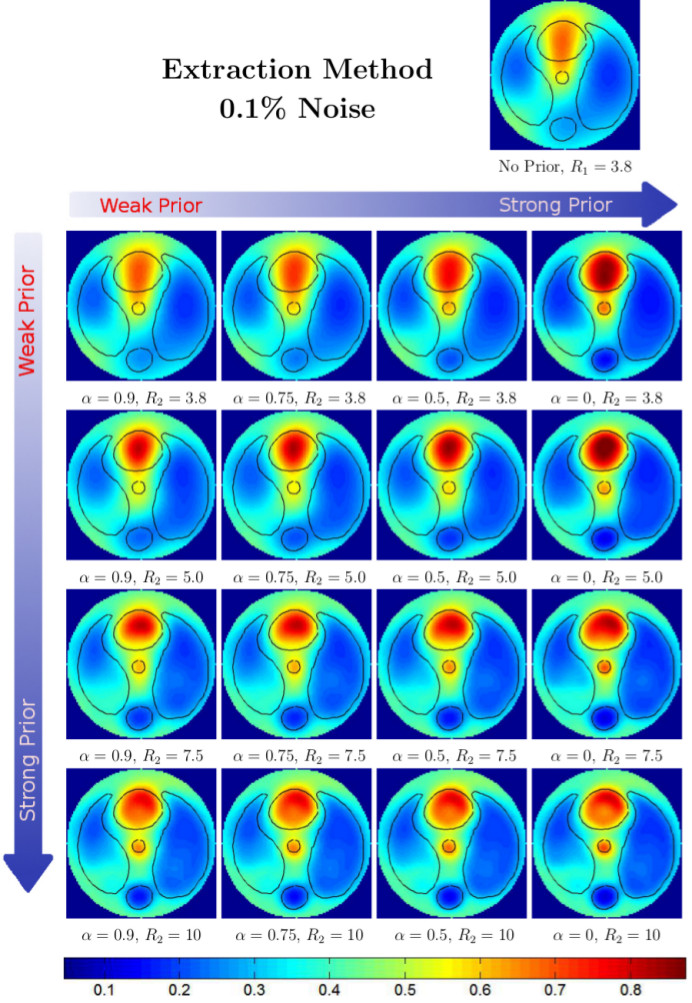}
  \caption{Reconstructions $\sigalpha$ for the 0.1\% noise case using the \emph{a priori} scheme with the extraction method, with various values of $\alpha$ and $R_2$. The reconstruction with no prior is at the top for comparison. The strength of the prior increases moving left to right and top to bottom. The scale bar at the bottom applies to all reconstructions at this noise level.}
\label{fig:Extraction_0p1}
\end{figure}

\begin{figure}[ht]
  \centering
    \includegraphics[height = 0.88\textheight]{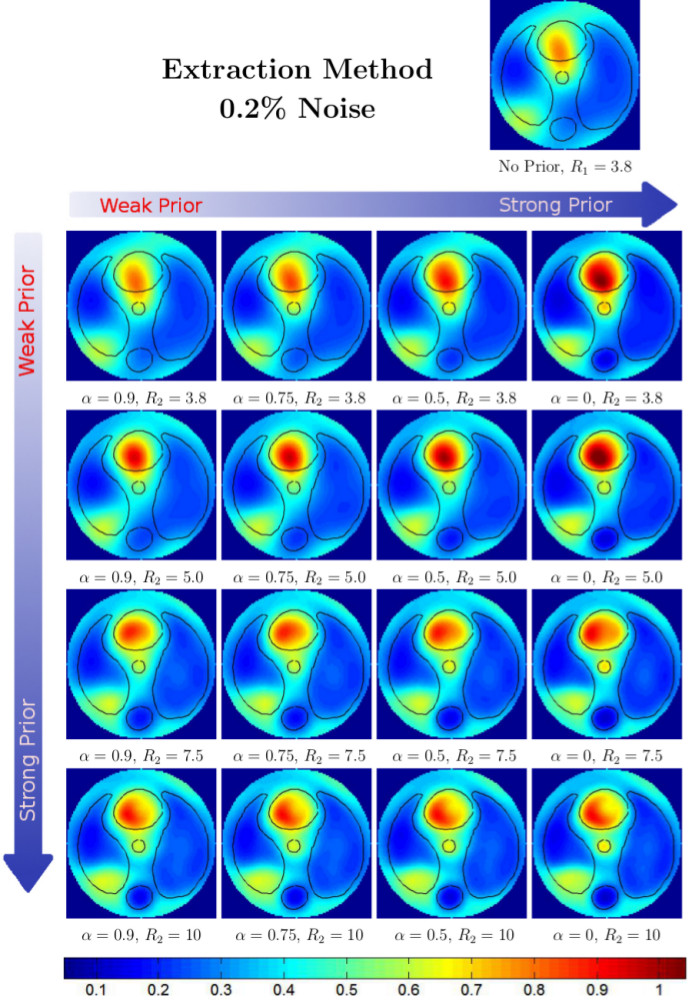}
  \caption{Reconstructions $\sigalpha$ for the 0.2\% noise case using the \emph{a priori} scheme with the extraction method, with various values of $\alpha$ and $R_2$. The reconstruction with no prior is at the top for comparison. The strength of the prior increases moving left to right and top to bottom. The scale bar at the bottom applies to all reconstructions at this noise level.}
\label{fig:Extraction_0p2}
\end{figure}

For the extraction of conductivity values, from the $\sigma$ reconstructions with the superimposed \emph{a priori} organ boundaries, it was clear that the left lung has increased conductivity toward the bottom, so we again divided the lung into top and bottom to construct the prior. For simplicity, we used the same dividing line as was used in \S\ref{BlindEstimateApplied}, and we selected the same values for $c$, $c_1$, and $c_2$ for all noise cases. Using the methods outlined in \S\ref{Extraction}, we extracted conductivity values from the reconstruction $\sigma$ to create $\sigpr$; these assigned values are given in Table~\ref{table:cond_ext}, along with the conductivity values used in the phantom for comparison. We then proceeded with the rest of the \emph{a priori} scheme outlined in \S\ref{APriori_outline}, using $R_1 = 3.8$, and testing various values for $R_2$ and $\alpha$. The resulting reconstructions are given in Figures~\ref{fig:Extraction_NoNoise}, \ref{fig:Extraction_0p1}, and \ref{fig:Extraction_0p2}.

In this method, with or without noise, the effusion is clearly visible for all weights of the prior, with improvement in the organ shapes as the weight of the prior increases.  Excellent reconstructions are found even in the presence of noise.  An iteration step is not included for this method since the first step produces very high quality reconstructions.

\section{Conclusion} \label{Conclusion}
A method for including {\em a priori} information in the 2-D D-bar algorithm was presented with two methods suggested for assigning conductivity values to the prior.  The {\em a priori} information is included in the scattering transform and in the integral equation for the CGO solution $\mu(z,k)$ and is weighted with two parameters $R_2$ and $\alpha$ in the scattering transform computation and the integral equation for $\mu$, respectively.  The method is demonstrated to be highly effective on numerically simulated data with noise levels typically used in EIT data simulations.  The method shows promise for clinical use in lung imaging when {\em a priori} information about organ boundaries can be obtained from a recent CT or MRI scan, for example.  Future work is needed to evaluate its clinical effectiveness.

\begin{table}[ht]
\centering
\caption{Conductivity values in S/m for the phantom as well as the ``blind estimate'' $\tsigpr$ values assigned, along with values assigned to $\tsigpr'$ in the subsequent iteration step, for each of the 3 noise cases.}
\label{table:conductivities2}
{\footnotesize
\begin{tabularx}{\textwidth}{|l|X|X|X|X|X|X|X|}
\hline
& Back- \newline ground & Heart & L Lung \newline top & L Lung \newline bottom &  R Lung & Aorta & Spine \\ \hline 
Used in phantom					& 0.424 & 0.750 & 0.240 &0.600  & 0.240 & 0.750  & 0.150  \\ \hline
``Blind estimates'' used in $\tsigpr$				& 0.500 & 0.800 & 0.200 &0.200 & 0.200 & 0.800  & 0.100  \\ \hline
Used in $\tsigpr'$, 0\% noise 			& 0.431 & 0.798 & 0.261 &0.364 & 0.233 & 0.798 & 0.187  \\ \hline
Used in $\tsigpr'$, 0.1\% noise 			& 0.427 & 0.767 & 0.274 &0.333 & 0.233 & 0.767 & 0.178  \\ \hline
Used in $\tsigpr'$, 0.2\% noise 			& 0.450 & 0.858 & 0.247 &0.428 & 0.247 & 0.858  & 0.163  \\ \hline
\end{tabularx}
}	
\end{table}

\begin{table}[ht]
\centering
\caption{Conductivity values in S/m for the phantom as well as $\tsigpr$ values computed using extraction method for each of the 3 noise cases.}
\label{table:cond_ext}
{\footnotesize
\begin{tabularx}{\textwidth}{|l|X|X|X|X|X|X|X|}
\hline
& Back- \newline ground & Heart & L Lung \newline top & L Lung \newline bottom &  R Lung & Aorta & Spine \\ \hline 
Used in phantom			& 0.424 & 0.750 & 0.240 &0.600  & 0.240 & 0.750  & 0.150  \\ \hline
Extracted from $\sigma$, 0\% noise 		& 0.401 & 0.681 & 0.283 &0.398  & 0.251 & 0.681  & 0.186  \\ \hline
Extracted from $\sigma$,  0.1\% noise 	& 0.393 & 0.648 & 0.292 &0.373  & 0.252 & 0.648 & 0.178  \\ \hline
Extracted from $\sigma$, 0.2\% noise 	& 0.423 & 0.742 & 0.272 &0.460  & 0.260 & 0.742  & 0.177  \\ \hline
\end{tabularx}
}	
\end{table}

\clearpage
\bibliographystyle{siam}
\bibliography{APrioriBib}
\end{document}